\documentclass[12pt]{amsart}
\usepackage{xcolor}
\usepackage{hyperref} 
\usepackage[utf8]{inputenc}
\usepackage[russian,english]{babel}
\usepackage{epsfig}
\usepackage{epigraph}
\author{Alexander Givental}
\title{On merit and equity in math}

\begin{document}
\begin{abstract} In recent paper {\em Quantifying Inequities and Documenting Elitism in PhD-granting Mathematical Sciences Departments in the United States} \cite{1}
by a group of accomplished and/or aspiring mathematicians, the authors use data to challenge ``the idea that the mathematical sciences in the United States is a meritocracy''. We show that, regardless of the validity of the challenged idea, the arguments in the paper are invalid. Namely, among other flaws, they rely on a logical trick of asserting the validity of the conclusion derived from an assumption which the authors neglect to test (and which actually contradicts their own data). \end{abstract}

\maketitle

My first refutation of paper \cite{1} is contained in another article \cite{2} (rejected by the arXiv, apparently, for its not entirely academic character), which uses a number of very elementary mathematical models to dispel several popular ways of ``social justice'' reasoning. The models are so elementary that I wouldn't bother writing a separate article about them, but the aim of exposing the fallacies of \cite{1} turned out to be a sufficiently worthy occasion. The present text is an abbreviated (and more ``academic'') version of \cite{2}, focusing on the latter task only.

The inequities, allegedly quantified in \cite{1}, concern: (a) underrepresentation of women among math faculty, and (b) disproportional over-funding of ``elite''math departments by the NSF's Department of Mathematical Sciences (DMS) program. Below we expose a number of elementary mathematical, logical, and factual errors committed in \cite{1}, which invalidate the methodology of the paper and undermine its conclusions. 

\medskip

{\bf 1.} According to the authors of \cite{1}, in the past decade, women account for $30\%$ of all new PhD degrees awarded in math, and therefore (they argue) in pure meritocracy the portion of women on the math faculty of PhD-granting mathematics department had to be also $30\%$. The banal fallacy of this reasoning can be described by a preschool-level math theorem: {\em the maximum is equal the average only when the quantity is constant}. Namely, since the portion of women awarded PhD in math has been {\em growing}, even if women form $30\%$ of new hires, women still wouldn't make 30\% of all math faculty including those hired in the past. For instance, in my home department (see \cite{3}) women currently make $13\%$ of all Senate faculty (which includes emeriti), $20\%$ of currently active, and $30\%$ among assistant and associate professors.                

\medskip

{\bf 2.} Even assuming (for the sake of the argument) that the quantity is constant, i.e. that women receive $30\%$ of PhD degrees in math for the span of an entire generation, the conclusion in the paper that ``math in the US is not a meritocracy'' would be based on explicitly formulated assumptions (see page 2 in \cite{1}) which can be summarized this way: It can be expected that women form a certain fraction (e.g. $30\%$) of all people interested in pursuing a particular career (e.g. in math), but then in pure meritocracy they would constitute that same fraction in every tier of the achievement ladder (such as PhD degree or  --- a higher tier --- professorship).

The obvious fallacy of conclusions derived from an assumption is that the assumption might be false. And this particular assumption is well-known to be generally false. Indeed, while it is possible that non-proportional attrition rates of men and women up the achievement ladder are due to biases, there are many examples where it would be hard to suspect any non-meritorious forces in play. A famous one is chess \cite{4} where among FIDE-rated players 11\% are women though among grandmasters they make $<3\%$. Another (random) example is found in a study \cite{5} of Germany's several-tier system of selection of high-schoolers to the International Chemical Olympiad, where the portion of girls decreased with each tier (in 2017 from 48\% to 0), and a similar decrease is found among the recipients of awards (honorable mention, bronze, silver, gold) at the Olympiads themselves.

Actually the very fact that women are being awarded $30\%$ of PhD degrees in math (rather than $50\%$) already invalidates the authors' assumption --- unless they also attribute this disproportion to bias against women (and hence bias against men in other subjects, because \cite{6} overall women in the US earn more PhD degrees than men each year starting 2006).  

We refer to \cite{2} for the discussion of one illustrative model showing that the hypothesis about the natural proportionality of the distributions of men and women over all tiers of the achievement ladder has no grounds whatsoever. Even under the assumption of equal standard deviations (true for distributions of FIDE-rated players \cite{4}, and useful anyway since the opposite costed one Harvard president their job \cite{7}), the distributions may have different means (as it happens, e.g with ACT math scores \cite{8}). In the example of normal distributions this leads to exponentially decaying ratio.

\medskip

{\bf 3.} Concerning DMS funding of mathematics departments, the authors of \cite{1} arrive (on page 9) at the conclusion that the funding ``is not meritocratic and reinforces elitism''. They base this conclusion on their finding that the top $25\%$ of the departments ordered by their {\em prestige ranking} (whatever this means) absorb $65\%$ of DMS funds allocated to all departments so ranked, and that the top $20\%$ of all awardees (PIs? departments? --- unclear) absorb $86\%$ of all DMS funds.

The obvious disconnect between the finding and conclusion manifests in the absence of any attempt to take into account the scientific and educational impact of those elite departments --- or merely their size! What if those most prestigious $25\%$ account for $65\%$ of all math faculty? Even if not, what if they account for the overwhelming majority of math PhD students (to be supported by the DMS funds)? What about their scientific impact (for evaluation of which there are any number of metrics, based on citation indexes, journal rankings, academy memberships, prizes and medals, etc. and --- yes --- NSF awards)? The authors speculate along the lines of ``the rich-get-richer'', i.e. that better funding is not a result of success in math but a primary cause of it. Yet they have absolutely no data supporting this fantasy.  

\medskip

{\bf 4.} Eventually, on page 10, in the spirit of {\em 1984}, they advocate for what amounts to imposing racial, ethnic, and gender quotas in hiring: ``we advocate for the redefinition of 'prestige' and 'elitism' in mathematics in a way that better reflects equity for excellence. For example, institutions that might deserve our respect are institutions that reflect the diversity of the US population among their faculty and doctoral graduates in mathematics.''

To illustrate this new $raison\ d'\hat{e}tre$ of mathematics departments, they present the list of most prestigious 10: {\em ``The top ten PhD-granting Mathematics departments by the percentage of women.''} In fact, as it is discussed with some detail in \cite{2}, while four of the ten can be estimated to yield altogether about 8 math PhDs per year (for comparison: UC Berkeley Math department awards 25-30), in the remaining six I didn't manage to locate any PhD programs in math at all.  

\bigskip

The ``non-academic'' question somehow addressed in \cite{2} is why this scientifically looking yet grotesquely misleading paper was written at all. It is inconceivable that a group of mathematicians could inadvertently make errors as naive as some noted above. Alternatively, one could insist that this is a typical case of cynical exploitation of the now-popular ideology for personal gain. Yet, in \cite{2}, exercising the presumption of innocence, I argue that the paper is a hoax, intended to compromise the DEI axiomatics, and to illuminate the sheer incompatibility of {\em meritocracy} with {\em equity}.

\enddocument